\documentstyle{amsppt}
\magnification=\magstep1
\parindent=1em
\baselineskip 15pt
\hsize=12.3 cm
\vsize=18.5 cm
\NoRunningHeads
\NoBlackBoxes
\pageno=1

%


\def\ker{\text{ker }}
\topmatter

\title
 Frames of translates
\endtitle
\author
Peter G. Casazza, Ole Christensen, and Nigel J. Kalton
\endauthor
\address1
Department of Mathematics,
The University of Missouri,
Columbia, Missouri 65211,
USA
\endaddress
\address2
Mathematical Institute,
Building 303,
Technical University of Denmark,
2800 Lyngby,
DENMARK
\endaddress
\address3
Department of Mathematics,
The University of Missouri,
Columbia, Missouri 65211,
USA
\endaddress
\email
pete\@casazza.math.missouri.edu; olechr\@mat.dtu.dk;
nigel\@math.missouri.edu
\endemail

\thanks
The first author was supported by NSF grant DMS 970618; the second author
by the Danish Research Council,
and the third author was
supported by
NSF grants DMS 9500125 and 9870027.  The second author would also like
to thank the University of Missouri-Columbia for its hospitality.
\endthanks
%
\abstract
Frames consisting of translates of a single function play an
important role
in sampling theory as well as in wavelet theory and Gabor analysis.
We give a
necessary and sufficient condition for a subfamily of regularly
spaced translates
of a function ${\phi}\in L^{2}(R)$, $({\tau}_{nb}{\phi})_{n\in
\Lambda}$, ${\Lambda}\subset Z$, to form a frame (resp. Riesz
basis)
for its closed linear span.  One consequence is that if ${\Lambda}
\subset
N$, then this family is a frame if and only if it is a Riesz basis.
 For the
case of arbitrary translates of a function ${\phi}\in L^{1}(R)$ we show
that for sparse sets, having an upper frame bound is equivalent to
the family being a frame sequence. Also we give some relatively
mild density conditions will
yield frame sequences.  Finally, we us the fractional
Hausdorff dimension to identify classes of exact frame sequences.
\endabstract

\endtopmatter
\document
\baselineskip=15pt

\heading{1. Preliminaries}
\endheading
\vskip10pt

Let $H$ be a Hilbert space with inner product $< \cdot, \cdot >$ linear
in the first entry. Let $\Lambda$ be any countable index set. Recall that
a sequence
$(f_n)_{n\in\Lambda}
\subseteq H$ is a
{\it frame for $H$} if
$$
\exists A,B>0: \ A ||f||^2 \le \sum_{n\in\Lambda} |<f,f_n>|^2 \le B
||f||^2
,
\
\forall f \in H. \tag1.1
$$
$A,B$ are called the {\it lower and upper frame bounds}.  They are not
unique:
the biggest lower bound and the smallest upper bound are called the
{\it optimal frame bounds}.  It can be shown that
every element $f \in H$ has a
representation as an infinite linear combination of the frame elements:
using the {\it frame operator}
$$S : H \to H, \ \ Sf= \sum <f,f_n>f_n.$$
Then
$$f= SS^{-1}f= \sum <f, S^{-1}f_n>f_n , \ \ \forall f \in H.$$
Clearly a frame $(f_n)$ is total, i.e., $\overline{span}
(f_n)=H.$ In case $(f_n)$ is a frame for $\overline{span}(f_n)$ (which in
general might be a proper subspace of $H$) we say that $(f_n)$ is a
{\it frame sequence}.  We say that $(f_{n})$ is an {\it exact frame
sequence} if $(f_{n})$ is a Riesz basis for $\overline{span}\
(f_{n})$; i.e., if
$$
\exists A,B>0:\ A\sum |c_{n}|^{2} \le \|\sum_{n\in\Lambda}
c_{n}f_{n}\|^{2}
\le B
\sum |c_{n}|^{2}, \tag 1.2
$$
for all $(c_{n})\in c_{00}(\Lambda)$, the space of finitely nonzero
sequences.

Let us remark here that (1.1) implies that $(f_n)_{n\in\Lambda}$ is a
frame for $H$ with frame bounds $A$ and $B$, if and only if there is
bounded operator
$V:H\to\ell^2(\Lambda)$ defined by $Vf=\{\langle
f,f_n\rangle\}_{n\in\Lambda}$  with $A\|f\|^2\le  \|Vf\|^2 \le B\|f\|^2$
for $f\in H.$  Thus $V$ must an isomorphism onto its range $V(H).$
Let $(e_n)_{n\in\Lambda}$ denote the canonical basis vectors in
$\ell^2(\lambda).$  We then see, taking adjoints, that the above
condition is also
equivalent to the existence  of a bounded surjective operator
$T:\ell_2(\Lambda)\to H$ so that $Te_n=f_n$ and
$A\|a\|^2 \le \|Ta\|^2\le B\|a\|^2$ for $a\in (\ker T)^{\perp}$ or,
equivalently $Ad(a,\ker T)^2\le \|Ta\|^2\le B\|f\|^2.$  The frame
operator is $S=TV=TT^*.$

Thus in particular $(f_n)_{n\in\Lambda}$ is a frame sequence with
constants $A,B$ if and only if we have the linear map
$T:c_{00}(\Lambda)\to H$ defined by $Te_n=f_n$ extends to a bounded
operator on $\ell_2(\Lambda)$ and satisfies
$$ A\|a\|^2\le \|Ta\|^2\le B\|a\|^2 \tag 1.3$$ for $a\in (\ker
T)^{\perp}.$  Comparing with (1.2) we see that $(f_n)$ is an exact frame
sequence if and only if $T$ is one-one.

The most important examples of frame sequences comes from sampling
theory \cite{5}.
For $x \in\Bbb  R$ we define {\it translation by $x$} by
$$
{\tau}_x : L^2(\Bbb R) \to L^2(\Bbb R), ({\tau}_xf)(y)=f(y-x), \ y \in
\Bbb R.
$$
It is proved in \cite{4} that a collection $( {\tau}_{x_n}{\phi})$,
where
${\phi} \in L^2(R), (x_n) \subseteq R$ can never be a frame for
$L^2(\Bbb R)$.
However, frame sequences of this form exists and play an important
role in
sampling theory as well as wavelet theory, cite{5,6,8}.

We start by considering sequences of the form
$({\tau}_{nb}{\phi})_{n\in {\Lambda}}$, ${\phi} \in L^{2}(R)$,
${\Lambda}\subset Z$,   In section 2
we give necessary and sufficient conditions for such sets to form frame
(resp. exact frame) sequences.  This extends work of Benedetto and
Walnut \cite{2} and Benedetto and Li \cite{1} as well as removing an
unnecessary
hypothesis in their results.  Kim and Lim \cite{12} also showed that this
was an unnecessary hypothesis in \cite{2}.  As one of several
applications of this result, we show that if ${\Lambda}\subset N$,
then $({\tau}_{nb}{\phi})_{n\in {\Lambda}}$ is a frame sequence if
and only if it is an exact frame sequence.
In section 3 we give conditions for an arbitrary sequence of
translates to
have finite upper frame bound.  In section 4 we consider  the
case where ${\Phi}_{b}$ (see section 2 for the definition) is continuous
and use the cardinality of the zero set of ${\Phi}_{b}$ and the
density of our set to produce frame sequences.  Finally,
in section 5 we relate the fractional Hausdorff dimension to exact frame
sequences of translates.

\heading{2.  Frames of Translates}
\endheading
\vskip10pt

Our first theorem is a generalization of a
result of Benedetto amd Li \cite{1}. The proof we give is considerably
simpler.

We first introduce some notation.  For a function $\phi\in L^1(\Bbb R)$
we denote by $\hat \phi$ the {\it Fourier transform of $\phi$}
$$ \hat \phi(\xi) = \int\phi(x) e^{-2\pi i\xi x}(x)dx.$$
As usual the definition of the Fourier transform extends to an isometry
$\phi\to\hat \phi$ on $L^2(\Bbb R).$

Now suppose $\phi\in L^2(\Bbb R)$ and that $b> 0.$ Let us identify
the
circle $\Bbb T$ with the interval $[0,1)$ via the standard map $\xi\to
e^{2\pi i\xi}.$ We define the function
$\Phi_b:\Bbb T\to \Bbb R$ by
$$ \Phi_b(\xi) =\sum_{n\in\Bbb Z}|\hat\phi(\frac{\xi+n}{b})|^2.$$
Note that $\Phi_b\in L^1(\Bbb T).$

For any $n\in\Bbb Z$ we note that
$$ \langle \tau_{nb}\phi,\phi\rangle
=\langle e^{-2\pi in\xi b}\hat\phi,\hat\phi\rangle
=\frac1b\int_0^1\Phi_b(\xi)e^{-2\pi in\xi}d\xi =\frac1b\hat\Phi_b(n).$$

If $\Lambda\subset\Bbb Z$ we let $H_{\Lambda}$ be the closed subspace of
$L^2(\Bbb T)$ generated by the characters $e^{2\pi in\xi}$ for $n\in
\Lambda.$  We let $E_{\Lambda}$ be the closed subspace of $H_{\Lambda}$
of all $f$ such that $\Phi_b(\xi)f(\xi)=0$ a.e.   If $f\in H_{\Lambda}$
we denote by $d(f,E_{\Lambda})$ the distance of $f$ to the subspace
$E_{\Lambda}.$

\proclaim{Theorem 2.1}Suppose $\phi\in L^2(\Bbb R)$ and $b> 0.$  If
$\Lambda\subset\Bbb Z$ then $(\tau_{nb}\phi)_{n\in\Lambda}$ is a frame
sequence with frame bounds $A$ and $B$ if and only if for every
 $f\in H_{\Lambda}$ we have
$$ Ad(f,E_{\Lambda})^2\le \frac1b \int_0^1|f(\xi)|^2\Phi_b(\xi)d\xi\le  B
\|f\|^2,$$ or equivalently, for all $f\in H_{\Lambda}\cap
E_{\Lambda}^{\perp}$,
$$ A\|f\|^2\le \frac1b \int_0^1|f(\xi)|^2\Phi_b(\xi)d\xi\le  B \|f\|^2.$$
Furthermore, if this condition is satisfied,
$(\tau_{nb}\phi)_{n\in\Lambda}$ is an exact frame sequence
with the same bounds if and only if $E_{\Lambda}=\{0\}.$\endproclaim

\demo{Proof}By our remarks and the definition
$(\tau_{nb}\phi)_{n\in\Lambda}$
is a frame sequence with frame bounds $A,B$ if and only if the linear map
$T:c_{00}(\Lambda)\to L^2(\Bbb R)$ defined by
$Te_n=\tau_{nb}\phi$
extends to a
 bounded linear operator $T:\ell^2(\Lambda)\to L^2(\Bbb R)$ such that
$Te_n=\tau_{nb}\phi$  and
$$ Ad(u,\ker T)^2  \le \|Tu\|^2\le B\|u\|^2$$ for $u\in \ell^2(\Lambda).$

Let $U:H_{\Lambda}\to\ell^2(\Lambda)$ be the natural isometry $Uf=\{\hat
f(n)\}_{n\in \Lambda}.$  Then for any trigonometric polynomial $f\in
H_{\Lambda}$
$$\align \|TUf\|^2 &=\|\sum_{n\in\Lambda}\hat f(n)\tau_{nb}\phi\|^2 \\
 &= \int_{-\infty}^{\infty}|\sum_{n\in\Lambda}\hat f(n)e^{-2\pi
inb\xi}\hat\phi(\xi)|^2d\xi\\
&=\frac1b \sum_{n\in\Bbb Z}\int_0^1 |f(\xi)|^2
|\hat\phi(\frac{n+\xi}{b})|^2d\xi\\
&= \frac1b\int_0^1|f(\xi)|^2\Phi_b(\xi)d\xi.\endalign$$

This immediately implies the theorem.\qed\enddemo

 Theorem 2.1 yields a generalization of
a result of Benedetto and Li \cite{1} which is part (3) of the next
theorem (note that an unnecesary hypothesis in \cite{1} is also
eliminated).

\proclaim{Theorem 2.2}
If $\phi \in L^{2}(\Bbb R)$, and $b>0$ then:

(1)  $({\tau}_{nb}{\phi})_{n\in \Bbb Z}$ is an orthonormal sequence if
and only if
$$
{\Phi}_{b}(\gamma) = b\ \ \text{a.e.}
$$

(2)  $({\tau}_{nb}{\phi})_{n\in \Bbb Z}$ is an exact frame sequence with
frame bounds $A,\ B$ if and
only if
$$
bA\le {\Phi}_{b}(\gamma)\le bB\ \ \text{a.e.}
$$

(3)  $({\tau}_{nb}{\phi})_{n\in \Bbb Z}$ is a frame sequence with frame
bounds $A,\
B$ if and only if
$$
bA\le {\Phi}_{b}(\gamma)  \le bB \ \ \text{a.e.}  \tag 2.1
$$
on $\Bbb T\setminus N_{b}$ where $N_{b} = \{\xi\in \Bbb T:
{\Phi}_{b}(\xi)
= 0\}$.
\endproclaim

\demo{Proof} Note that (1) follows easily from the fact that
$(\tau_{nb}\phi)_{n\in \Bbb Z}$ is orthonormal if and only if $TU$ is
unitary.  (2) is immediate from (3). For (3) we note that if
$\Lambda=\Bbb Z$ then
$H_{\Lambda}=L^2(\Bbb T)$ and
$E_{\Lambda}=L^2(N_b)$.  Hence $d(f,E_{\Lambda})^2=\int_{\Bbb T\setminus
N_b}|f|^2d\xi.$  The proof is then immediate.\qed\enddemo

The next theorem addresses the relationship between frame
properties for two
families $({\tau}_{nb} \phi)_{n\in \Bbb Z}$ and $({\tau}_{na}
\phi)_{n\in \Bbb Z}$.
Observe that we can assume $0<b<a$ without loss of generality.

\proclaim{Theorem 2.3}
Let $a,b\in \Bbb R$ with $0 < b < a$.

(1)  There exists a function $\phi \in L^{2}(\Bbb R)$ so that
$({\tau}_{na} \phi)_{n\in \Bbb Z}$ is a frame sequence but
$({\tau}_{nb} \phi)_{n\in \Bbb Z}$ is not a frame sequence.

(2)  If $\frac{a}{b}\in \Bbb Z$ and $({\tau}_{nb} \phi)_{n\in \Bbb Z}$ is a frame
sequence, then
$({\tau}_{na} \phi)_{n\in \Bbb Z}$ is a frame sequence.

(3)  If $\frac{a}{b}\notin \Bbb Z$ then there is a function $\phi
\in L^{2}(\Bbb R)$ so that $({\tau}_{nb} \phi)_{n\in \Bbb Z}$ is a frame
sequence, but
$({\tau}_{na} \phi)_{n\in \Bbb Z}$ is not a frame sequence.

\endproclaim

\demo{Proof}
(1)  Define a function $\phi \in L^{2}(R)$ so that
$$
\hat{\phi}(\xi) = 1,\ \ \text{if}\ \ \ 0\le \xi \le \frac{1}{a},
$$
$$
\hat{\phi}(\xi) = \frac{ab}{b-a}(x-\frac{1}{b})
  \ \ \text{if}\ \ \ \frac{1}{a}\le \xi \le \frac{1}{b},
$$
and
$$
\hat{\phi}(\xi) = 0,\ \ \text{otherwise}.
$$
Then we can easily check that
$$
{\Phi}_{b}(\xi) =
 \hat{\phi}(\frac{\xi}{b})^{2},
$$
for all $\xi \in \Bbb T$.  Hence, ${\Phi}_{b}(\xi)$ is not
bounded below on
$\Bbb T$ and so $({\tau}_{nb} \phi)_{n\in \Bbb Z}$ is not a frame sequence.  On the
other hand,
if $\xi \in \Bbb T$ then $\frac{\xi}{a}\in [0,\frac{1}{a}]$, so
that ${\Phi}_{a}(\xi)
\ge 1$.  Therefore, ${\Phi}_{a}(\xi)$ is bounded below on $\Bbb T$
and so
$({\tau}_{na}\phi)_{n\in \Bbb Z}$ is a frame sequence.

(2)  By our assumption, there is a natural number $m$ so that $a = mb$.
Hence,
$$
{\Phi}_{a}(\xi) = \sum_{n}|\hat{\phi}(\frac{\xi + n}{a})|^{2} =
\sum_{n}|\hat{\phi}(\frac{\xi + n}{mb})|^{2}
$$
$$
= \sum_{k=0}^{m-1}|\hat{\phi}(\frac{\frac{\xi}{m} + \frac{k}{m}
+ n}{b})|^{2} = \sum_{k=0}^{m-1}{\Phi}_{b}(\frac{\xi + k}{m}).
$$
It follows that if $({\tau}_{nb} \phi)_{n\in \Bbb Z}$ is a frame
sequence, then either
${\Phi}_{a}(\xi)= 0$ or $A\le {\Phi}_{a}(\xi)\le Bm$.  Hence,
$({\tau}_{na} \phi)_{n\in \Bbb Z}$ is also a frame sequence.

(3)  Since $\frac{a}{b}\notin \Bbb Z$, there is an $0 < \epsilon $ so that
$$
\frac{\xi +n}{a} \notin [\frac{1}{b},\frac{1}{b}+\epsilon],\ \
\forall
n\in \Bbb Z,\ \ \forall 0<\xi < \epsilon.
$$
Define a function $\phi \in L^{2}(\Bbb R)$ so that,
$$
\hat{\phi}(\xi) = \xi, \ \ \forall \ \ 0\le \xi \le
\frac{1}{a},
$$
$$
\hat{\phi}(\xi) = 1, \ \ \forall \ \ \xi \in
 [\frac{1}{b},\frac{1}{b}+\epsilon],
$$
and $\hat{\phi}(\xi) = 0$ otherwise.  Then,
$$
{\Phi}_{a}(\xi) = \hat{\phi}(\frac{\xi}{a})^{2},\ \ \forall
\ \ 0\le \gamma
\le {\epsilon},
$$
and so ${\Phi}_{a}(\xi)$ is not bounded below on $\Bbb T$.
Therefore,
 $({\tau}_{na} \phi)_{n\in \Bbb Z}$ is not a frame sequence.  On the
other hand,
$$
{\Phi}_{b}(\xi) \ge 1\ \  \text{on} \ \ [0,c],
$$
where $c = \text{min}\{\epsilon,\frac{\epsilon}{b}\}$,
and for $\xi \in [c,1]$, either ${\Phi}_{b}(\xi) = 0$
or ${\Phi}_{b}(\xi) \ge {c}^{2}$.
It follows that where ${\Phi}_{b}(\xi)$ is non-zero, it is
bounded above
and below and hence $({\tau}_{nb} \phi)_{n\in \Bbb Z}$ is a frame
sequence.  \qed

\enddemo

Notice that Theorem 2.3(2) implies that if $(\tau_{nb}\phi)_{n\in\Bbb Z}$
is a frame sequence than $(\tau_{nb}\phi)_{n\in\Lambda}$ is a frame
sequence whenever $\Lambda$ is a {\it subgroup} of $\Bbb Z.$
It is also clear that any subsequence of an exact frame sequence remains
an exact frame sequence.  The following result is a converse of this.

\proclaim{Theorem 2.4}Suppose $\phi\in L^2(\Bbb R)$ and $b>0.$  Then the
sequence $(\tau_{nb}\phi)_{n=1}^{\infty}$  is a frame sequence if and
only if $(\tau_{nb}\phi)_{n\in\Bbb Z}$ is an exact frame
sequence.\endproclaim

\demo{Proof}Assume $(\tau_{nb}\phi)_{n\in\Bbb N}$ is a frame sequence.
Then (see for example \cite{7}) if $0\neq f\in H_{\Bbb N}$ we have
that $\log|f|\in L^1$ and so, in particular, $|f|>0$ a.e.  This implies
that $E_{\Bbb N}=\{0\}.$  It follows that if $A,B$ are the frame bounds
for $(\tau_{nb})_{n\in \Bbb N}$ then
for every trigometric polynomial in $H_{\Bbb N}$ we have
$$A\|f\|^2 \le \int |f(\xi)|^2\Phi_b(\xi)d\xi\le B\|f\|^2.$$

Now suppose $f$ is any trigonometric polynomial.  Then for large enough
$n$ we have that $e^{2\pi ni\xi}f\in H_{\Bbb N}$.
Thus the same inequality follows trivially for all trigonometric
polynomials in $L^2(\Bbb T).$   This implies the theorem.\qed\enddemo

\demo{Remark} If $(\tau_{nb}\phi)_{n\in \Bbb Z}$ is a frame sequence but
not an
exact sequence frame, then of course $(\tau_{nb}\phi)_{n=1}^{\infty}$
cannot be a
frame sequence.  The set $(\tau_{nb}\phi)_{n\in \Bbb Z}$ is
linearly independent, but it follows
easily that the lower frame bound for $(\tau_{nb}\phi)_{n=1}^N$ must
converge to $0$ as $N\to\infty.$  \enddemo

\demo{Remark} Note that if $\Lambda\subset \Bbb N$ then the argument of
the above theorem shows that $(\tau_{nb}\phi)_{n\in\Lambda}$ is a frame
sequence if and only if it is also an exact frame sequence.\enddemo

\heading{3.  Upper frame bounds for subsequences
}
\endheading
\vskip10pt

In this section we give some criteria for the existence of an upper frame
bound for a sequence $(\tau_{nb}\phi)_{n\in \Lambda}.$

First we introduce some notation.  If $\Lambda$ is a countable subset of
$\Bbb R$
  we define for $x>0$ the function
$$ D_{\Lambda}(x)=\sup_{t\in\Bbb R}|\Lambda \cap [t,t+x]|.$$

Our first result concerns general conditions for an upper frame bound to
exist for a sequence of translates (which need not in this case be
regularly spaced).

\proclaim{Theorem 3.1}Let $F:(0,\infty)\to(0,\infty)$ be a
monotone-decreasing function such that $F\in L_2.$  Define
$$ G(x)=F(x)\int_0^xF(t)dt + \int_x^{\infty}F(t)^2dt.\tag 3.1$$
Let $\Lambda=(\lambda_n)$ be a countable subset of $\Bbb R$.
Then:
\newline
(1) Suppose $\phi\in L^2(\Bbb R)$ is such that for some $X$ and every $x$
with $|x|\ge X$ we have
$|\phi(x)|\le F(|x|).$
Then in order that $(\tau_{\lambda_n}\phi)$ have an upper frame bound it
is sufficient that
$$ \int_1^{\infty}G(x)D_{\Lambda}(x)\frac{dx}{x} <\infty.$$
\newline
(2) If $\psi(x)=F(|x|)$ for $x\neq 0$ then a necessary condition for
$(\tau_{\lambda_n}\psi)$ to have an upper frame bound is that
$$ \sup_{x>1}G(x)D_{\Lambda}(x)<\infty.$$
\endproclaim

\demo{Remark} Note that $G$ is continuous, decreasing and bounded by
$G(0)=\lim_{x\to 0}G(x)=\int_0^{\infty}F(x)^2dx.$\enddemo

\demo{Proof}Suppose $t>0.$  Note that
$$ \int_{-\infty}^{\infty}\psi(x+t)\psi(x-t)dx=
2\int_{0}^{2t}\psi(x+t)\psi(x-t)dx +
2\int_{2t}^{\infty}\psi(x+t)\psi(x-t)dx.$$
Now
$$2F(3t)\int_0^tF(x)dx \le \int_0^{2t}\psi(x+t)\psi(x-t)dx \le
2F(t)\int_0^t F(x) dx$$
and
$$
\int_{3t}^{\infty}F(x)^2dx \le \int_{2t}^{\infty}\psi(x+t)\psi(x-t)dx
\le
\int_{t}^{\infty}F(x)^2dx. $$

 We now argue that
$$ \frac{4}{3}G(3t) \le
 \int_{-\infty}^{\infty}\psi(x+t)\psi(x-t)dx \le 4G(t).$$
Indeed
we have
$$
\align
\frac43 G(3t) &= \frac43 F(3t)\int_0^{3t}F(x)dx
+\frac43\int_{3t}^{\infty}F(x)^2dx\\
&\le 4F(3t)\int_0^tF(x)dx
+\frac43\int_{2t}^{\infty}\psi(x+t)\psi(x-t)dx\\
&\le 2\int_0^{2t}\psi(x+t)\psi(x-t)dx +
\frac43\int_{2t}^{\infty}\psi(x+t)\psi(x-t)dx\\
& \le \int_{-\infty}^{\infty}\psi(x+t)\psi(x-t)dx\\
&= 2\int_0^{2t}\psi(x+t)\psi(x-t)dx
+2\int_{2t}^{\infty}\psi(x+t)\psi(x-t)dx\\
&\le 4F(t)\int_0^tF(x)dx + 2\int_t^{\infty}F(x)^2dx \\
&\le 4G(t).\endalign$$

Let us now prove (1). We first estimate $|\phi|\le
|\psi|+|\phi\chi_{[-X,X]}|.$  It now follows that if $t\ge X$
$$ \int_{-\infty}^{\infty} |\phi(x+t)||\phi(x-t)|dx \le
4G(t)+CF(2t-X)$$ for a suitable constant $C$.
It follows (noting that $G$ is bounded) that for a suitable constant $C$
we have that
$$ \langle \phi,\tau_a\phi\rangle \le CG(\frac12a)$$
for all $a>0.$

 Consider any finitely nonzero sequence $(c_n)$.
Then
$$ \align
\|\sum_{n=1}^{\infty}c_n\tau_{\lambda_n}\phi\|^2
&=\sum_{m=1}^{\infty}\sum_{n=1}^{\infty}\overline{c_m}c_n\langle
\phi,\tau_{\lambda_n-\lambda_m}\phi\rangle\\
&\le
C\sum_{m=1}^{\infty}\sum_{n=1}^{\infty}|c_m||c_n|G(\frac12|\lambda_m-
\lambda_n|)\\
&\le
\frac12 C\sum_{m=1}^{\infty}\sum_{n=1}^{\infty}(|c_m|^2+|c_n|^2)
G(\frac12|\lambda_m-\lambda_n|)\\
&\le C
\sum_{m=1}^{\infty}|c_m|^2\sum_{n=1}^{\infty}G(\frac12|\lambda_m-
\lambda_n|).\endalign $$
Now
$$\align \sum_{n=1}^{\infty}G(\frac12|\lambda_m-\lambda_n|) &\le  2
\sum_{k=0}^{\infty}\sum_{2^{k-1}<|\lambda_m-\lambda_n|\le
2^k}G(2^{k-1})+ \sum_{|\lambda_m-\lambda_n|\le 1/2}G(0)\\
&\le
4\sum_{k=0}^{\infty}D_{\Lambda}(2^{k-1})G(2^{k-1})+D_{\Lambda}(1)G(0)\\
&\le
8\sum_{k=0}^{\infty}D_{\Lambda}(2^{k-2})G(2^{k-1})+D_{\lambda}(1)G(0)\\
&\le 8(\log 2)^{-1}\int_{1/4}^{\infty}D_{\lambda}(x)G(x)\frac{dx}{x} +
D_{\Lambda}(1)G(0)\\
&\le 8(\log 2)^{-1}\int_1^{\infty}D_{\Lambda}(x)G(x)\frac{dx}{x}+ (1+\log
4)D_{\Lambda}(1)G(0).\endalign $$
This establishes (1).

To prove (2) fix $x>1$ and $t\in\Bbb R,$ and let $\Bbb A=\{n:t\le
\lambda_n\le t+x\}.$   Then
$$ \align
\|\sum_{n\in\Bbb A}\tau_{\lambda_n}\phi\|^2 &\ge\frac{4}{3}\sum_{m\in\Bbb
A}\sum_{n\in
\Bbb A} G(\frac{3}{2}|\lambda_m-\lambda_n|)\\
&\ge \frac{4}{3}|\Bbb A|^2 G(\frac{3}{2}x)\endalign $$
so that if $B$ is the upper frame bound,
$$ G(\frac{3}{2}x)|\Bbb A|\le \frac{3}{4}B.$$
Hence
$$ G(\frac{3}{2}x)D_{\Lambda}(x) \le \frac{3}{4}B.$$
Now this implies that
$$ G(x) D_{\Lambda}(x) \le 2G(x)D_{\Lambda}(\frac12 x)\le \frac32 B.$$
This completes the proof.\qed\enddemo

\proclaim{Corollary 3.2}Under the hypotheses of the theorem if $F\in
L^1(\Bbb R)$ then a necessary and sufficient condition for
$(\tau_{\lambda_n}\phi)$ to have an upper frame bound is that
 $D_{\Lambda}(1)<\infty.$
\endproclaim
\demo{Proof}In this case we clearly have an estimate  $G(x)\le CF(x)$.
The condition
$D_{\Lambda}(1)<\infty$ is equivalent to $D_{\Lambda}(x)\le Cx$ for $x\ge
1.$\enddemo

\proclaim{Lemma 3.3}Suppose that $F:(0,\infty)\to(0,\infty)$ satisfies
the conditions of Theorem 3.1 and is such that for some $\epsilon>0$ we
have that
$x^{1-\epsilon}F(x)$ is increasing for $x\ge 1$ and
$x^{1+\epsilon}F(x)^2$
is decreasing for $x\ge 1.$  Then there is a constant $C$ so that
$C^{-1}xF(x)^2
\le G(x)\le CxF(x)^2,$ for $x\ge1$, where $G$ is defined in
(3.1).\endproclaim

\demo{Proof} This is an immediate calculation from (3.1).\enddemo

\demo{Remark} Note that in the theorem if $F(x)=\min(1,x^{-a})$ where
$\frac12<a<1$ then $G(x)\approx \min(1,x^{1-2a})$.  Hence a
necessary condition
for $(\tau_{\lambda_n}(\phi))$ to have an upper-frame bound is that
$D_{\Lambda}(x)\le Cx^{2a-1}$ for $x\ge 1$ and a sufficient condition
is that $D_{\Lambda}(x)\le Cx^{2a-1-\epsilon}$ for $x\ge 1$, for some
$\epsilon>0.$
If $F(x)=\min(1,x^{-1})$ then $G(x)\approx x^{-1}(1+\log x)$ for $x\ge
1$.\enddemo

We now prove a more precise result for the case when $F$ is of the form
given in Lemma 3.3, by giving a condition which is close to necessary
and sufficient. This result will not be used later and the reader who
is only interested in later results may therefore omit it. The
argument is simply an adaptation of an argument
used for a similar result in \cite{9} and derives from more
general results in \cite{10}.

\proclaim{Theorem 3.4}Suppose $F:(0,\infty)\to(0,\infty)$ is a monotone
decreasing function with the property that for some $\epsilon>0$ we have
$x^{1-\epsilon}F(x)$ is increasing for $x\ge 1$ but $x^{1+\epsilon}F(x)$
is decreasing for $x\ge 1$.  Let $\Lambda=(\lambda_n)$ be a countable
subset of $\Bbb R.$  Then:\newline
(1) Suppose $\phi\in L^2(\Bbb R)$ is such that for some $X$ and every
$x\ge X$ we have $|\phi(x)|\le F(|x|).$  Then in order that
$(\tau_{\lambda_n}\phi)$ has an upper frame bound it is necessary that
there is a constant $C$ so that for every finite interval $I$ we have
 $$ \sum_{\lambda_m,\lambda_n\in I}G(|\lambda_m-\lambda_n|) \le
C|\Lambda\cap [t,t+x]| \tag 3.2$$
where $G$ is defined by (3.1).
\newline
(2) If $\psi(x)=F(|x|)$ then (3.2) is also necessary for
$(\tau_{\lambda_n}\phi)$ to have an upper frame-bound.
\endproclaim

\demo{Proof}
 In this proof, we will use
$C$ for a constant depending only
on $F$ which may vary from line to line.  Note that we have estimates of
the form $F(x/2)\le CF(x)$ and $G(x/2)\le CG(x).$  It follows therefore
that
$$ \frac1C G(|a|)\le \langle \psi,\tau_a\phi\rangle\le CG(|a|).$$
Let us first prove (2) (necessity).  Indeed for any bounded
interval $I,$ if
 $A=\{n:\lambda_n\in I\},$  then
$$ \|\sum_{n\in A}\tau_{\lambda_n}\psi\|_2^2 \ge \frac1C \sum_{m,n\in
A}G(|\lambda_m-\lambda_n|)$$
and so the conclusion is immediate.

The other direction (1) is harder.
It will be convenient to assume $F(x)=1$ for $0\le x\le 1.$
Clearly if we prove the result under this assumption the general case
will follow.

Let us start by observing the estimate
from (3.2) that if $A=\{n:\lambda_n\in I\}$ then
$ G(x) |A|^2 \le C|A|$ so that $$|\Lambda\cap I|\le
CG(x)^{-1}\tag3.3.$$

Now fix $t\in \Bbb R$ and suppose $x>1$.  We estimate, using Lemma 3.3:
$$ \align\sum_{|t-\lambda_n|\ge x}F(|t-\lambda_n|)&\le
C\sum_{k=1}^{\infty}F(2^{k-1}x)|\{n:2^{k-1}x\le |t-\lambda_n|\le 2^kx\}|
\\
&\le C\sum_{k=0}^{\infty}F(2^kx)G(2^kx)^{-1}\\
&\le C \sum_{k=0}^{\infty}2^{-k}x^{-1}F(2^kx)^{-1}\\
&\le C\sum_{k=0}^{\infty}2^{-\epsilon
k}(2^{{\epsilon-1}k}x^{-1}F(2^kx)^{-1})\\
&\le Cx^{-1}F(x)^{-1}.
\endalign
$$
In general for all $x>0$ we conclude an estimate of the type:
$$ \sum_{|t-\lambda_n|\ge x}F(|t-\lambda_n|)\le
C\min(1,x^{-1})F(x)^{-1}.\tag3.4$$

We now introduce two further functions.  Let $N(t)=|\Lambda\cap
[t-1,t+1]|$ and let $H(t)=\sum_{n}F(|t-\lambda_n|).$  Our assumptions
give us that $N(t)\le C$ and $N(t)\le H(t).$  By (3.4) we
have
$H(t)\le N(t)+C.$

Suppose we have an interval $I=[x-h/2,x+h/2]$ of length $h>1.$ Let
$J=[x-h,x+h]$. Then
we can write $H(t)=H_1(t)+H_2(t)$ where $H_1(t)=\sum_{\lambda_n\in
J
}F(|t-\lambda_n|)$ and $H_2(t)=H(t)-H_1(t).$  Using (3.4) we have, if
$|t-x|\le \frac12h,$
$$\align H_2(t)^2 &\le 2\sum_{\lambda_m\notin J}\sum_{|\lambda_n-t|\ge
|\lambda_m-t|}F(|t-\lambda_m|)F(|t-\lambda_n|) \\
&\le C\sum_{\lambda_m\notin J}|t-\lambda_m|^{-1}.\endalign $$
On the other hand:
$$ \int_{-\infty}^{\infty} H_1(t)^2  \le C\sum_{\lambda_m,\lambda_n\in
J}G(|\lambda_m-\lambda_n|) \le C\sum_{\lambda_m\in J}1.$$
Combining we have
$$ \int_IH(t)^2dt \le C\sum_{m}\min(1,h|x-\lambda_m|^{-1}).$$
We use this inequality to estimate
$$\align \int_{-\infty}^{\infty}F(|t-x|)H(x)^2dx &\le
\int_{t-1}^{t+1}H(x)^2dx +
C\sum_{k=1}^{\infty}F(2^k)\int_{2^{k-1}<|t-x|<2^k}H(x)^2dx\\
&\le C\sum_{k=0}^{\infty}F(2^k)\sum_m \min(1,2^k|t-\lambda_m|^{-1})\\
&= C\sum_m \sum_{k=0}^{\infty}F(2^k)\min(1,2^k|t-\lambda_m|^{-1})\\
& \le
C\sum_m\sum_{2^k<|x-\lambda_m|}2^kF(2^k)|t-\lambda_m|+C\sum_m\sum_{2^k\ge
|x-\lambda_m|}F(2^k)\\
&\le C\sum_m F(|t-\lambda_m|)=CH(t) \tag 3.5 \endalign $$

For each $n$ let $E_n=(\lambda_n-\frac12,\lambda_n+\frac12).$  For any
finitely nonzero sequence $(a_n)$, let $f=\sum_n |a_n|\chi_{E_n}$. Then
since $f^2\le N(t) \sum_n|a_n|^2\chi_{E_n}$ we have $\|f\|^2\le
C(\sum_n|a_n|^2).$  We also have:
$$\align \|\sum_na_n\tau_{\lambda_n}\phi\|_2^2 &\le
C\sum_{m,n}|a_m||a_n|\int_{-\infty}^{\infty}F(|t-\lambda_m|)F(|t-\lambda_n|)dt\\
&\le
C\sum_{m,n}|a_m||a_n|\int_{-\infty}^{\infty}\int_{E_n}F(|t-x|)dx
\int_{E_m} F(|t-y|)dy\,dt \\
&\le C\int_{-\infty}^{\infty}(\tilde f(t))^2dt, \tag 3.6\endalign $$
where $$ \tilde f(t)=\int_{-\infty}^{\infty}F(|t-s|)f(s)ds.$$

Now suppose $g$ is any nonnegative $L^2-$function of compact support with
$\|g\|=1$ and suppose $\tilde
g(t)=\int_{-\infty}^{\infty}F(|t-s|)g(s)ds.$ Then
$\tilde g\in L_2$ and
$$
\align
\int_{\infty}^{\infty}g(t)\tilde f(t)dt &=\int_{-\infty}^{\infty}\tilde
g(t)f(t)dt\\
&\le C\int_{-\infty}^{\infty}\tilde g(t)f(t)H(t)^2dt\\
&\le C(\int_{-\infty}^{\infty}\tilde g(t)^2H(t)^2dt)^{1/2}
(\sum_n|a_n|^2)^{1/2}.\tag 3.7
\endalign $$

Now
$$\align \tilde g(t)^2
&=\int_{-\infty}^{\infty}\int_{-\infty}^{\infty}F(|t-x|)F(|t-y|)g(x)g(y)dx\,dy\\
&\le
C\int_{-\infty}^{\infty}\int_{-\infty}^{\infty}F(|y-x|)(F(|t-x|)+F(|t-y|))g(x)g(y)dx\,dy\\
&\le
C\int_{-\infty}^{\infty}\int_{-\infty}^{\infty}F(|y-x|)F(|t-x|)g(x)g(y)dx\,dy\\
&\le
C\int_{-\infty}^{\infty}F(|t-x|)g(x)\tilde g(x)dx.\tag 3.8\endalign $$
Hence, using (3.5),
$$ \align
\int_{-\infty}^{\infty}\tilde g(t)^2H(t)dt &\le
C\int_{-\infty}^{\infty}\int_{-\infty}^{\infty}F(|t-x|)g(x)\tilde
g(x)H(t)^2 dt\,dx
\\
&\le \int_{-\infty}^{\infty}g(x)\tilde g(x)
\int_{-\infty}^{\infty}F(|x-t|)H(t)^2dt\\
&\le C\int_{-\infty}^{\infty}g(x)\tilde g(x) H(x)dx\\
&\le C\left(\int_{-\infty}^{\infty}\tilde g(x)^2H(x)^2dx\right)^{1/2}
\endalign $$
which implies
$$ \int_{-\infty}^{\infty}\tilde g(x)^2 H(x)^2 dx \le C.$$
Now (3.7) implies
$$ \int_{-\infty}^{\infty}g(t)\tilde f(t)dt \le C(\sum_m|a_m|^2)^{1/2}$$
and as $g$ is arbitrary (subject to being of compact support and norm
one) this implies that
$\|\tilde f\|^2 \le C\sum_m|a_m|^2$ and (3.6) then yields the
theorem.\qed\enddemo

For regularly spaced sequences we have some other criteria. We use the
terminology of Section 2.  Recall that a subset $\Lambda$ of $\Bbb Z$ is
a $\Lambda(p)-$set (for $p>2$) if the $L_2-$ and $L_p-$norms are
equivalent on $H_{\Lambda}.$ See \cite{13} and \cite{3}

\proclaim{Theorem 3.5}Suppose $b\in\Bbb R$ and that $\phi\in L^2(\Bbb
R).$  Suppose $\Lambda$ is a $\Lambda(p)-$set where $p>2$ and that
$\Phi_b\in L_r(\Bbb T)$ where $\frac1r+\frac2p=1$.  Then
$(\tau_{nb}\phi)_{n\in\Lambda}$ has an upper-frame bound.
\endproclaim

\demo{Proof}If $f\in H_{\Lambda}$, we observe
that
$$ \int_{\Bbb T}|f(\gamma)|^2\Phi_b(\gamma)d\gamma \le \|f\|^2_p
\|\Phi_b\|_r$$ and so the theorem is immediate.\qed\enddemo

\demo{Remark}In order that $\Phi_b\in L_r$ it is sufficient by the
Hausdorff-Young inequality that $\sum_{n\in\Bbb Z}|\langle
\phi,\tau_{nb}\phi\rangle|^{p/2}<\infty$.\enddemo

\heading{4. Frame sequences when $\phi\in L^1(\Bbb R).$}
\endheading
\vskip10pt

Let us call a subset $\Lambda$ of $\Bbb Z$ {\it sparse} if for any
$n\in \Bbb N$ the set $\Lambda \cap (\Lambda +n)$ is finite.  Thus any
increasing sequence $(\lambda_n)$ is sparse if and only if
$\lim_{n\to\infty}(\lambda_n-\lambda_{n-1})=\infty.$

\proclaim{Theorem 4.1}Suppose $0\neq \phi\in L^1(\Bbb R)\cap L^2(\Bbb
R)$ and $\Lambda$ is a sparse subset of $\Bbb Z$.
Then
 $(\tau_{nb}\phi)_{n\in\Lambda}$ has an upper frame bound if and only if
it is a frame.  In particular if $|\phi(x)|=O(F(x))$ where
$F:(0,\infty)\to(0,\infty)$ is monotone-decreasing and integrable then
$(\tau_{nb}\phi)_{n\in\Lambda}$ is a frame.
 \endproclaim

\demo{Proof}Suppose $(\tau_{nb}\phi)_{n\in\Lambda}$ has an upper frame
bound but no lower frame bound.  We appeal to Theorem 2.1.  There is a
sequence $(f_n)\in H_{\Lambda}\cap E_{\Lambda}^{\perp}$ so that
$$\lim_{n\to\infty}\int_0^1|f_n
(\xi)|^2\Phi_b(\xi)d\xi=0$$ but $\|f_n\|=1$
for all
$n.$  By passing to a subsequence, we can assume without loss of
generality that
$f_n$ converges weakly to some $f\in H_{\Lambda}\cap
E_{\Lambda}^{\perp},$ and even further that
$g_n=\frac1n(f_1+\cdots+f_n)$ converges in norm to $f$.
Since $(g_n)$ converges to $f$ in measure, it follows from Fatou's Lemma
that $\int_0^1|f(\xi)|^2\Phi_b(\xi)d\xi=0$ and so $f\in E_{\Lambda}.$
Hence $f=0,$ and $f_n$ converges weakly to zero.

Now we estimate
$$ |\int_0^1 |f_n(\xi)|^2e^{-2\pi i k\xi}d\xi|\le\sum_{j\in\Bbb Z}|\hat
f_n(j)||\hat f_n(k+j)|.$$
Hence
$$ |\int_0^1 |f_n(\xi)|^2e^{-2\pi i k\xi}d\xi|\le\sum_{j\in
\Lambda\cap (\Lambda-k)}|\hat f_n(j)||\hat f_n(j+k)|.$$
Since $f_n$ converges to $0$ weakly this converges to $0$ if $k\neq 0.$
Hence the measures $|f_n(\xi)|^2d\xi$ converge weak$^*$ to
$d\xi$ in $C(\Bbb T)^*$.

Now since $\phi\in L^1$ the function $\hat \phi$ is continuous and so
$\Phi_b(t)$ is lower-semi-continuous and $2\pi-$periodic when regarded as
a function on $\Bbb R.$  In particular $\Phi_b$ is lower-semi-continuous
on $\Bbb T$ and hence there is a sequence of functions $\psi_n\in C(\Bbb
T)$ such that $0\le \psi_n\uparrow \Phi_b$ pointwise.

Clearly
$$\align \int_0^1\psi_k(\xi)d\xi
&=\lim_{n\to\infty}\int_0^1\psi_k(\xi)|f_n(\xi)|^2d\xi
\\
&\le \limsup_{n\to\infty}\int_0^1\Phi_b(\xi)|f_n(\xi)|^2d\xi\\
&=0.\endalign $$
Hence $\int_0^1\Phi_b(\xi)d\xi=0$ which is a contradiction.
\qed\enddemo

If $\phi$ decays rapidly enough at $\infty$ we can achieve a stronger
result. We will need the following lemma.

\proclaim{Lemma 4.2}Suppose $\Lambda\subset \Bbb Z$,
 and suppose $f\in H_{\Lambda}$, with $\|f\|=1.$ Let
$F=|f|^2.$
 If $\Bbb J$ is a finite interval in $\Bbb Z$ then
$$ \sum_{n\in\Bbb J}|\hat F(n)| \le  D_{\Lambda}(|\Bbb J|)$$
\endproclaim
\demo{Proof} We have
$$ |\hat F(n)|\le \sum_{k\in\Bbb
Z}\chi_{\Lambda}(k)\chi_{\Lambda}(n-k)|\hat f(k)||\hat f(n-k)|.$$
Making the estimate $|\hat f(k)||\hat f(n-k)|\le \frac12 (|\hat f(n)|^2
+|\hat f(n-k)|^2$ we obtain
$$ |\hat F(n)| \le \sum_{k\in \Bbb Z}|\hat f(k)|^2 \chi_{\Lambda}(n-k)$$
so that
$$ \sum_{n\in\Bbb J}|\hat F(n)| \le D_{\Lambda}(|\Bbb J|)\sum_{k\in\Bbb
Z} |\hat f(k)|^2$$ which immediately yields the lemma.\qed\enddemo

\proclaim{Theorem 4.3}Suppose $\phi\in L^2(\Bbb R)$ and that $b\in\Bbb
R.$  Assume $\Phi_b$ is continuous on $\Bbb T$ and has exactly $N$ zeros.
Suppose $\Lambda$ is any set satisfying the density condition
$$ \lim_{x\to\infty}\frac{D_{\Lambda}(x)}{x}<\frac1N.$$  Then
$(\tau_{nb}\phi)_{n\in\Lambda}$ is a frame sequence.
\endproclaim

\demo{Proof}Since $\Phi_b$ is bounded the upper frame bound is trivial.
Assume there is no lower frame bound.  Then there is a sequence $f_n\in
 H_{\Lambda}$ with $\|f_n\|=1$ and $$\lim_{n\to\infty}\int \Phi_b(x)
|f(x)|^2 =0.$$
Without loss of generality we can suppose the measures $|f_n(x)|^2dx $
converge weak$^*$ to a probability measure $d\mu$.  Since $\Phi_b$ is
continuous (lower semi-continuity would suffice!) we have
$$ \int \Phi_b(x)d\mu(x)=0$$ so that $\mu$ can be written in the form
$\mu=\sum_{k=1}^Na_k\delta_{t_k}$ where $t_1,\ldots,t_N$ are the zeros of
$\Phi_b$ and $0\le a_k$ with $\sum a_k=1.$

Now if $F_n=|f_n|^2$ we have, applying Lemma 4.2, for every natural
number
$m$ that
$$ \sum_{k=-m}^m|\hat F_n(k)|^2 \le D_{\Lambda}(2m+1).$$
Hence
$$ \sum_{k=-m}^m|\hat \mu(k)|^2 \le D_{\Lambda}(2m+1).$$
It follows from a theorem of Wiener \cite{11} that
$$ \sum_{k=1}^Na_k^2 \le
\lim_{m\to\infty}\frac{D_{\Lambda}(2m+1)}{2m+1}<\frac1N.$$
But then this contradicts the Cauchy-Schwartz inequality since
$\sum_{k=1}^Na_k=1.$\qed\enddemo

\proclaim{Theorem 4.4}Suppose $h:(0,\infty)\to (0,\infty)$ is an
increasing  continuous function such that there exist constants
$1<c_1<c_2<\infty$ with $c_1 h(x)\le h(2x)\le c_2h(x)$ for all $x,$ and
$\int_1^{\infty}t^{-2}h(t)dt=\infty.$
Suppose $\phi\in L_2(\Bbb R)$ satisfies the
condition that $\phi(x)=O(e^{-\delta h(x)}))$ as $x\to\infty$, where
$\delta>0.$
 Then, for
any $b\in\Bbb R,$ there exists an integer $N$ so that if $\Lambda\subset
\Bbb Z$ with
$\lim_{x\to\infty}\frac{D_{\Lambda}(x)}{x}<\frac1N$ then
$(\tau_{nb}\phi)_{n\in\Lambda}$ is a frame sequence.
\endproclaim

\demo{Proof}The idea of the proof is to show that the function $\Phi_b$
is $C^{\infty}$ on $\Bbb T$ and is {\it quasi-analytic} (see
\cite{11})
so that there is no $\xi\in [0,1)$ for which all derivatives
$\Phi_b^{(n)}(\xi)=0$ for $n\ge 0.$  It then follows from Taylor's
theorem that the zeros of $\Phi_b$ are isolated and hence finite.
Theorem 4.3 will then complete the proof.  To show $\Phi_b$ is
quasi-analytic we need to show that if $M_n=\|\Phi_b^{(n)}\|_2$ and
$\tau(r)=\inf_n(M_nr^n)$ then
$$ \int_0^{\infty}\frac{\log \tau(r)}{1+r^2}dr=-\infty.\tag4.1$$

Throughout the proof $C$ will denote a constant, depending on
$h$ and $\phi$, which may vary from line to line but is independent of
$x,m,n,k$ etc.
We begin with the observation that there exists a constant $C$ so that
$$ \frac1C\int_0^x h(t)\frac{dt}{t} \le h(x) \le
C\int_0^xh(t)\frac{dt}{t}$$
so that by replacing $h$ with $H(x)=\int_0^xh(t)\frac{dt}{t}$ we can
suppose $h$ is continuously differentiable and that there is a constant
$C$ so that $$\frac1C\le \frac{xh'(x)}{h(x)}\le C$$ for all $x>0.$
If we let $v=h^{-1}$ be the inverse function then we also have
$$\frac1C\le \frac{xv'(x)}{v(x)}\le C.$$

Let $m\in\Bbb N$.  Then $x^me^{-h(x)}$ attains a maximum at a point $x_m$
where
$x_mh'(x_m)=m$ and so $C^{-1}m\le h(x_m)\le Cm.$  Hence also
$C^{-1}v(m)\le x_m\le Cv(m).$  Combining we have
$$ x^me^{-h(x)}\le (Cv(m))^m.$$
Let $N$ be an integer greater than $2/\delta.$ It follows that
$$ x^{Nm}e^{-h(x)}\le (C(v(Nm)))^{Nm}$$
and so
$$ x^m e^{-\delta h(x)/2} \le (Cv(m))^m.$$
>From this we obtain
$$ \int_0^{\infty}x^me^{-\delta h(x)}dx \le (Cv(m))^m
\int_0^{\infty}e^{-\delta h(x)/2}dx \le (Cv(m))^m.$$

We now use the argument of Theorem 3.1 to deduce that
$$|\langle\phi,\tau_{nb}\phi\rangle| \le C\min
(1,e^{-\delta h(\frac12|n|b)})$$ for some suitable constant $C$ and all
$n\in\Bbb Z.$  It follows easily that
$$\Phi_b(\xi)=\sum_{n\in\Bbb Z}\langle \phi,\tau_{nb}\phi\rangle e^{2\pi
in\xi}$$ is $C^{\infty}.$  Furthermore
$$ M_n=\|\Phi_b^{(n)}\| \le C\sum_{k\in\Bbb Z}|k|^n e^{-\delta
h(\frac12|k|b)}$$
for $n\in\Bbb N.$  Hence
$$ M_n\le C^n \int_0^{\infty}x^ne^{-\delta h(x)}dx \le (Cv(n))^n.$$
It follows that
$$ \tau(r) \le \inf_{n\ge 1} (Cv(n)r)^n.$$
For given $r<1$ choose $n=[Ch(r^{-1})/e].$ Then $Cv(n)r\le 1/e$ and
$n\ge Ch(r^{-1})/e-1.$  Hence
$$ \log\tau(r) \le -Ch(r^{-1})$$
for small enough $r.$  Now
$$
\int_0^{1}\frac{h(r^{-1})}{1+r^2}dr=\int_1^{\infty}\frac{h(r)}{1+r^2}dr
=\infty$$
so that (4.1) holds and the proof is complete.\qed\enddemo

\heading{5. Fractional dimension of the zero-set and frames}\endheading
\vskip10pt

For our next lemma (which is probably known) we let $\ell(I)$ be the
length of an interval $I\subset \Bbb T.$

\proclaim{Lemma 5.1}Suppose $\Lambda\subset \Bbb Z$,
 and suppose $f\in H_{\Lambda}$, with $\|f\|=1.$ Then
 there is a constant $C$ so that if $I$ is an interval contained in
$\Bbb T$ then
$$  \int_I |f(\xi)|^2 dx \le C\ell(I)D_{\Lambda}(\ell(I)^{-1}).$$
\endproclaim

\demo{Proof}
 Let $I$ be the interval $[t-h,t+h]$. Set $r=1-h$ and then let
$$ \psi(x) =\frac{1-r^2}{1-2r\cos(2\pi(x-t))+r^2}=\sum_{k\in \Bbb Z}
r^{|n|}e^{2\pi in(t-x)}.$$
Then for some absolute constant $C$ independent of $t,h$ we have
$\chi_I\le Ch^{-1}\psi.$   Hence
$$ \align\int_I |f(x)|^2 dx &\le Ch^{-1}\int_0^1 \psi(x)F(x) dx  \\
&\le Ch^{-1}\sum_{k\in \Bbb Z} r^{|k|}|\hat F(n-k)|\\
&\le 2Ch^{-1}\sum_{k=0}^{\infty}r^{mk}D_{\Lambda}(m) \endalign $$
for any integer $m.$  Let $m$ be chosen so that $\ell(I)^{-1}\le m\le
2\ell(I)^{-1}.$   Then $r^m\le (1-h)^{1/(2h)}\le c<1$ for some
absolute constant $c$.  We thus quickly obtain the result since
$D_{\Lambda}(2x)\le 2D_{\Lambda}(x)$ for any $x.$\qed\enddemo

Now if $0<\alpha<1$ let us define the {\it essential $\alpha$-Hausdorff
measure} denoted $H_{\alpha}(E)$ of a subset $E$ of $\Bbb T$  to be the
infimum of $\sum_{n=1}^{\infty}\ell(I_n)^{\alpha}$ over all collections
of intervals $(I_n)_{n=1}^{\infty}$ so that $E\subset F\cup
\bigcup_{n=1}^{\infty}I_n$ where $F$ has measure zero.

\proclaim{Theorem 5.2}Suppose $\phi\in L^2(\Bbb R)$, $b>0$ and
$\Lambda\subset \Bbb Z$ are such that $(\tau_{nb}\phi)_{n\in \Lambda}$
has an upper frame bound.  Suppose further that $0<\alpha<1$ and
$\lim_{\epsilon\to 0}H_{\alpha}(\Phi_b>\epsilon)=0$ and that
$D_{\Lambda}(x)\le Cx^{1-\alpha}$ for some $C$ and all $x.$   Then
$(\tau_{nb}\phi)_{n\in\Lambda}$ is an exact frame sequence. \endproclaim

\demo{Proof}If $\epsilon>0$ then pick intervals $I_n$ so that
$\{\Phi_b>\epsilon\}\subset \cup_{n=1}^{\infty}I_n$ up to a set of
measure zero and
$\sum_{n=1}^{\infty}\ell(I_n)^{\alpha}<2H_{\alpha}(\Phi_b>\epsilon).$
Then
$$ \int_0^1 \Phi_b(\xi)|f(\xi)|^2 \ge \epsilon(\|f\|^2
-\sum_{n=1}^{\infty}\int_{I_n}|f(\xi)|^2d\xi).$$
If $f\in H_{\Lambda}$ then
$$ \sum_{n=1}^{\infty}\int_{I_n}|f(\xi)|^2d\xi \le
C\|f\|^2\sum_{n=1}^{\infty}\ell(I_n)^{\alpha} \le 2
CH_{\alpha}(\Phi_B>\epsilon)\|f\|^2.$$
If we take $\epsilon$ small enough we have:
$$ \sum_{n=1}^{\infty}\int_{I_n}|f(\xi)|^2d\xi \le \frac12\|f\|^2$$
and the result follows by appealing to Theorem 2.1.\qed\enddemo

Combining this with Theorem 3.1 gives us the following theorem:

\proclaim{Theorem 5.3}Suppose $\frac12<a<1$ and that $\phi\in L^2(\Bbb
R)$ satisfies the conditions that $|\phi(x)|=O(x^{-a-\epsilon})$ as
$|x|\to \infty$ where $\epsilon>0$ and that $\lim_{\epsilon\to
0}H_{2a-1}(\Phi_b>\epsilon)=0.$   Then for any subset $\Lambda\subset
\Bbb Z$ with $D_{\Lambda}(x) \le Cx^{2(1-a)}$ we have that
$(\tau_{nb}\phi)_{n\in\Lambda}$ is an exact frame sequence.\endproclaim

Note that in the case when $\Phi_b$ is lower-semi-continuous the
condition $\lim_{\epsilon\to 0}H_{\alpha}(\Phi_b>\epsilon)$ reduces to
the condition that $h_{\alpha}(\Phi_b=0)=0$ where $h_{\alpha}(E)$ is the
infimum of $\sum\ell(I_n)^{\alpha}$ over all coverings $E\subset
\cup_{n=1}^{\infty}I_n.$  Thus in this case the condition is essentially
a condition on the fractal dimension of the zero set of $\Phi_b.$

\demo{Example}
We conclude by constructing an example where $\Phi_b$ is bounded,
$\lim_{\epsilon\to 0}H_{\alpha}(\Phi_b>\epsilon)=0,$
for every $\epsilon>0$ we have an estimate
$D_{\Lambda}(x)\le C_{\epsilon}x^{1-\alpha+\epsilon},$  but
$(\tau_{nb}\phi)_{n\in\Lambda}$ fails to be a frame sequence.
We observe
however we do not have any such example with $\Phi_b$ continuous or
lower-semicontinuous.

For each $n\in\Bbb N$ let $m_n=\max([\alpha n-\sqrt n],0).$  We
then let
$\Lambda =\cup_{n\ge 1}\{2^{n}+k2^{m_n}: 1\le k\le 2^{n-m_n}\}.$
Now suppose $N=2^p$.  Then for any interval $\Bbb J$ of length $N$ it is
clear that we have
$$ |\Lambda\cap J| \le \sum_{j=1}^{p-1}2^{j-m_j} + 2|\Bbb J|\max_{j\ge
p}2^{-m_j}.$$  Hence
$$ |\Lambda\cap J| \le 2^{p-m_p}\sum_{j=1}^{p-1}2^{-j+m_p-m_{p-j}} +
2^{p-m_p}\max_{j\ge p}2^{m_p-m_j}.$$
Since we have an estimate $m_{p}\le m_{p-j}+C\beta j$ where
$\alpha<\beta<1$ and $C$ is a constant independent of $j,p,$ this implies
that
$$ D_{\Lambda}(2^p) \le C2^{p-m_p}$$
which implies $D_{\Lambda}(x)\le C_{\epsilon}x^{1-\alpha+\epsilon}$ for
all $\epsilon>0.$

Now suppose
$$f_n(\xi)=2^{(m_n-n)/2}\sum_{k=1}^{2^{n-m_n}}e^{2\pi ik2^{m_n}\xi}$$ so
that
$\|f_n\|=1$ and $f_n\in H_{\Lambda}.$  Notice that
$$ |f_n(\xi)| = 2^{(m_n-n)/2} \left|\frac{\sin
(2^{n-1}\xi/2)}{\sin(2^{m_n-1}\xi)}\right|.$$

  Let
$E_n=\{|f_n|<2^{\frac12(n-m_n-\frac12\sqrt n)}\},$ and let $F_n=\Bbb
T\setminus E_n.$ We will
have for an appropriate constant $C$ that $$\int_{E_n}|f(\xi)|^2d\xi \le
C2^{-\frac12\sqrt n}.$$

Also
$F_n$ is a union
of at most $2^{m_n}$ equal intervals and has total measure bounded by
$C2^{m_n-n+\frac12\sqrt n}.$  Hence
$$H_{\alpha}(F_n) \le C2^{m_n-\alpha n+\frac12\alpha\sqrt n)} \le
2^{-\frac12\alpha\sqrt n}.$$

Set $F_0=\Bbb T$. Now define $\Phi$ on $\Bbb T=[0,1)$ by
$$ \Phi(\xi)=\inf_{k\ge 0}2^{-k}\chi_{F_k}.$$
Then $H_{\alpha}(\cup_{j\ge k}F_j)\to 0$ so that $\Phi>0$ a.e. and indeed
$\lim_{\epsilon>0}H_{\alpha}(\Phi>\epsilon)=0.$

We choose $\phi\in L^2(\Bbb R)$ so that
$\hat\phi=\Phi^{1/2}\chi_{[0,1)}.$ If we take $b=1$ then
$\Phi=\Phi_b.$ Now
$$ \int_0^1|f_n(\xi)|^2\Phi(\xi)d\xi \le 2^{-n}\int_{F_n}|f(\xi)|^2d\xi +
C.2^{-\frac12\sqrt n} \to 0.$$
Hence $(\tau_n\phi)_{n\in\Lambda}$ cannot be a frame
sequence.\qed\enddemo

\enddocument